  \RequirePackage{lineno} 
\documentclass[12pt]{article}
\usepackage{amssymb}

\usepackage{amssymb,amsfonts}

\usepackage{pdflscape}
\usepackage{graphicx}
\usepackage{mathtools}
\usepackage{amsthm}
\usepackage{setspace}
\usepackage{diagbox}
\usepackage{adjustbox}
\usepackage{array}
\usepackage{multirow}
\usepackage{booktabs}
\usepackage[all,arc]{xy}
\usepackage{enumerate}
\usepackage{mathrsfs}
\usepackage{setspace}
\newtheorem{theorem}{Theorem}
\newtheorem{definition}{Definition}

\theoremstyle{definition}

\theoremstyle{remark}

\makeatletter
\let\c@equation\c@thm
\makeatother
\numberwithin{equation}{section}

\title{Formation of coalition structures as a  non-cooperative game 1: theory}

\author{Dmitry Levando, \thanks{ This paper comes as   a further development of the 3-rd and the 4-th chapters  of my PhD Thesis "Essays on Trade and Cooperation" at Ca Foscari University, Venezia, Italy. Acknowledge:  Nick Baigent, Phillip Bich, Jean-Marc Bonnisseau, Alex Boulatov, Emiliano Catonini, Giulio Codognato, Sergio Currarini, Luca Gelsomini,  Izhak Gilboa, Olga Gorelkina, Piero Gottardi, Roman Gorpenko, Eran Hanany, Mark Kelbert, Ludovic Julien, Alex Kokorev, Dmitry Makarov, Francois  Maniquet, Igal Miltaich, Stephane Menozzi,   Roger Myerson, Miklos Pinter, Ariel Rubinstein, Marina Sandomirskaya, William Thompson, Konstantin Sonin, Simone Tonin, Dimitrios Tsomocos, Eyal Winter, Shmuel Zamir. Special thanks to Fuad Aleskerov,   Shlomo Weber and Lev Gelman. Many thanks for  advices and for discussion to participants of  SI\&GE-2015, 2016 (an earlier version of the title was "A generalized Nash equilibrium"), CORE 50  Conference, CEPET 2016  Workshop, Games 2016 Congress, Games and Applications 2016 at Lisbon, Games and Optimization at St Etienne.  All possible mistakes are only mine. This research did not receive any specific grant from funding agencies in the public, commercial, or not-for-profit sectors. \newline E-mail for correspondence: dlevando (at)  hse.ru.}}
\date{} 

\begin{document}

%\thanks{ This paper comes as   a development of the 3-rd and the 4-th chapters  of my PhD Thesis "Essays on Trade and Cooperation" at Ca Foscari University, Venezia, Italy. Acknowledge:  Nick Baigent, Phillip Bich, Alex Boulatov, Emiliano Catonini, Giulio Codognato, Sergio Currarini, Luca Gelsomini,  Izhak Gilboa, Piero Gottardi, Roman Gorpenko, Eran Hanany, Ludovic Julien, Alex Kokorev, Dmitry Makarov, Francois  Maniquet, Igal Miltaich, Stephane Menozzi,   Roger Myerson, Ariel Rubinstein, Marina Sandomirskaya, Konstantin Sonin, Simone Tonin, Dimitrios Tsomocos, Eyal Winter, Shmuel Zamir. Special thanks to Fuad Aleskerov,   Shlomo Weber and Lev Gelman. Many thanks for  advices and for discussion to participants of  SIGE-2015 (an earlier version of the title was "A generalized Nash equilibrium"), CORE 50  Conference and CEPET 2016  Workshop and Games 2016 Congress.  All possible mistakes are only mine. %This research did not receive any specific grant from funding agencies in the public, commercial, or not-for-profit sectors. }
%\author{Dmitry V. Levando\thanks{Address: NRU-HSE, Moscow, Russian Federation, Shabolovka, 26, off. 214, e-mail: \textit{dlevando ( at )hse.ru}. }}
%\date{}                                           % Activate to display a given date or no date
\maketitle
\begin{abstract}
% ADVICE FROM TSOMOCOS
The paper defines a non-cooperative simultaneous finite game  to study coalition structure formation with intra and inter-coalition externalities.  The   novelty of the game is that the  game definition embeds a \textit{coalition structure formation mechanism}. 

This  mechanism portions a set of strategies of the game into partition-specific strategy domains, what makes every partition to be a non-cooperative game with partition-specific payoffs for every player.

The mechanism  includes a maximum coalition size, a set of eligible partitions with coalitions sizes no greater than this  number (which also serves as a restriction for a maximum number of deviators) and a coalition structure formation rule. 
The paper defines a family of nested  non-cooperative games parametrized by a size of a maximum coalition size. 

Every  game in the family has an equilibrium in mixed strategies. The equilibrium can  generate more than one coalition and encompasses intra and inter group externalities, what makes it different from the Shapley value. Presence of individual payoff allocation  makes it different from a strong Nash, coalition-proof equilibrium,  and some  other equilibrium concepts. The accompanying  papers demonstrate applications of the proposed toolkit.

 \end{abstract}

%\vspace{-35pt}

\noindent{\bf Keywords: } Noncooperative Games
%\smallskip
%\vspace{5pt}

%\begin{center}
%WORKING PAPER
%\\
%Not for distribution or circulation
%\end{center}
%\tableofcontents
%\end{document}
 
%\end{document}

%\bigskip

%\end{document}
%\pagebreak\doublespacing

%\section{}
%\subsection{}
\section{Subject of the paper}

The research topic of this paper was inspired by John Nash{'}s   \textquotedblleft Equilibrium Points in n-person games \textquotedblright  (1950). This remarkably short, but highly influential note of only  5 paragraphs  established an equilibrium concept and the proof of its existence  which did not require an explicit specification of a final coalition structure for a set of players. Prior to Nash{'}s paper, the generalization of the concept of equilibrium provided by von Neumann for the case of two-players zero-sum game was done by portioning the players into two groups and regarding several players as a single player.   However up to now, non of these approaches   resulted into expected  progress in    studying intra- and inter group externalities between players.  
%This paper studies the     question coming from  the paper  "Equilibrium points in n-person games"    by John Nash (1950). His paper established an equilibrium condition % for an assumption that   every player considers himself to be a unique deviator without an explicit specification of a final coalition structure for a set of players. 
%If we start thinking about multiple deviators then  we can not avoid coalition formation. 
%The  natural questions are:   if a player  assumes that beside her/him there are also some other deviators, then   which game is played,  and what  is an equilibrium concept for this game?

In this paper  "coalition structure", or  "partition"\footnote{Existing literature uses both terms.}  for short, is a collection of non-overlapping subsets from a set of players, which  in a union  make the original set. A group, or a coalition, is an element of a coalition structure or  of a  partition.\footnote{the same.}

A partition induces two types of effects on a player's payoff. The first, through actions of players of  the same coalition. These   effects  will be  addressed as \textit{intra}-coalition externalities.
The second, from all other players, who  are outside  coalition, and     belong to  different coalitions. This effects will be addressed as \textit{inter}-group externalities.

 Nash (1953)  suggested that cooperation   should be studied  within a group and in terms of non-cooperative fundamentals. This suggestion  is known now as the Nash Program.  Cooperative behavior  was understood as an activity inside  a group with positive externalities between players.  Nash did not  write   explicitly    about multi-coalition framework or coalition structures. Coalition structures allow us to study inter-coalition externalities,  along with intra-coalition, ones and to separate cooperation in payoffs and cooperation in allocation of players.
 
The best analogy for the difference between   Nash Program and   the current research is the  difference between partial and  general strategic equilibrium analysis\footnote{meaning, strategic market games of Shapley and Shubik} in economics. The former isolates a market  ignoring cross-market interactions, the latter explicitly studies cross market interactions.

The research agenda of the paper is:  how to construct coalition structures from  actions of self-interested agents. Moreover,
the paper offers  a  generalization of a non-cooperative game from Nash (1950) to address the problem  of coalition structure formation absent  in Nash (1953).
 The contributions of this paper are: a construction of a non-cooperative game with an  embedding  coalition structure formation mechanism, and   a parametrization of all constructed games by a number of deviators.  The three  subsequent  papers will demonstrate various  applications  of the model to different types of games, including construction of a non-cooperative stability criterium and defining a cooperation.  %,  how to measure stability of appearing coalition structures   in terms of a self-interested behavior and applications to network games.
 
 The paper has the following structure: Section 2 presents an example, on why studying inter-coalition externalities requires including coalition  structures into individual strategy sets, Section 3 presents a general model of the game, Section 4 presents an  example of the game. The last  section discusses the approach used and the results of the paper. 
\section{An example: corporate dinner game}
The simple example below shows why to study coalition structures we need to incorporate them into a strategy set, what is different from Hart and Kurz (1983). They  studied coalition formation in terms of  choice of a coalition by every player, not of a coalition structure.

Consider a game of 4 players: A  is  a President; B is a senior vice-president; $C_1,C_2$ are  two other vice-presidents. They have  a corporate dinner. \footnote{The game is designed  to obtain pure strategies equilibrium for some of  the players. Compare  this game with  a lunch game in the next paper, where all players are identical, and the result is a stochastic game, where states are equilibrium coalition structures.}
A coalition  is  a group of players   at \textit{one} table. Every player may sit only at one table.  A coalition structure is an  allocation of  all players  over no more than four tables. Empty tables are not   taken into account.

 Individual set of strategies is  a set of all coalition structures  for the   players, i.e. a set of all possible allocations of players over 4 non-empty tables. A player chooses one coalition structure, an element from his strategy set. % as different coalition structures have an impact on one{'}s payoff even if a coalition of the player is the same.  
  A set of  strategies in the game is a direct product of four individual strategy sets.  The choice of all players is a point in the set of strategies of the game. 
  
Preferences of the players  are  such that  everyone (besides A) would like   to have a dinner with A, but A only with B.  Everyone wants players outside his table to eat individually,  due to possible dissipation of rumors or information exchange. No one can   enforce  others to form or not to form coalitions.

  In  every coalition structure (or a partition) any   coalition (i.e. a table) is formed only if everybody at the table agrees to have dinner together, otherwise a non--selected player  eats alone. Further we address this as a coalition structure formation rule or a rule for simplicity. The same coalition  may belong  to different coalition structures, but with different allocations of players beyond it, compare lines 1 and 2 in Table \ref{dinner} further.  %A mismatch in choices  results in a dinner alone.%Family of rules divides a set of all strategies into coalition structures with fixed strategies, each one operating as a game with dummy strategies.
% A final realization of coalition structures depends on choices of  coalition structures of all players. 
 
  The game is simultaneous and one shot. A realization  of a final partition ( a coalition structure) depends on choices of coalition structures of  all players.  A choice of a player is an element of his/her strategy set.
  Example. Let player A choose $\{\{A,B \}, \{ C_1\},  \{ C_2\} \}$;  player  B choose $\{\{A,B \}, \{ C_1\},  \{ C_2\} \}$; player  $C_1$ choose $\{\{A,C_1 \}, \{ B\},  \{ C_2\} \}$,      and player  $C_2$ choose $\{\{A,C_2 \}, \{B\},  \{ C_2\} \}$. Then the final partition is  $\{\{A,B \}, \{ C_1\},  \{ C_2\} \}$.  It is clear that a strong Nash equilibrium (Aumann, 1960), which is based on  a deviation of a coalition of any size, does not  discriminate between the coalition structures mentioned above.
  
  Payoff profile of all players in the game should be defined for every final coalition structure. 
Table \ref{dinner} presents coalition structures only  with the best individual payoffs.\footnote{All other coalition structures have significantly  lower payoffs. } Thus only some partitions from the big set of all strategies deserves attention. The first column is  a number of a strategy. The second column is an allocation of players over coalition structures, and also is  a list of  the best\footnote{Payoffs in the rest coalition structures are much lower, and they are excluded from consideration.} final coalition structures from a set of all strategies. The third column is  an individual  payoff profile of all players. The fourth column is a list of values for coalitions in coalition structures if to calculate values using cooperative game theory.

  \begin{table}[htp]
\caption{Strategies and payoffs in the corporate dinner game} 
\begin{center}
\begin{tabular}{|c|c|c|c|}
\hline 
{\small{num}} & {\small{Best final  partitions}} & 
\begin{tabular}{c}  
{\small{ Non-cooperative payoff profile }} 
\\
{\small{ $(U_A,U_B,U_{C_1},U_{C_2})$ }} 
\end{tabular}
   & \begin{tabular}{c}  
{\small{   Values of coalitions as in  }} 
\\
{\small{  cooperative  game theory}} 
\end{tabular} \\
\hline
1&  $\{ \{ A,B\}, \{ C_1\}, \{ C_2\}  \}$ & (10,10,3,3)  & $20_{AB}$, $3_{C_1}$, $3_{C_2}$\\
 $2^*$& $\{ \{ A,B\}, \{ C_1, C_2\}  \}$ & (8,8,5,5)  & $16_{AB}$, $10_{C_1,C_2}$\\
 3&$\{ \{ A,C_1\}, \{ B, C_2\}  \}$ & (3,5,10,5)  & $13_{AC_1}$, $10_{BC_2}$ \\
  4&$\{ \{ A,C_1\}, \{ B\} , \{C_2\}  \}$ & (3,3,10,3) & $13_{AC_1}$, $3_B$, $3_{C_2}$ \\
 5&$\{ \{ A,C_2\}, \{ B, C_1\}  \}$ & (3,5,5,10)  &  $13_{AC_2}$, $10_{BC_1}$   \\
 6&  $\{ \{ A,C_2\}, \{ B\} , \{C_1\}  \}$ & (3,3,3,10) & $13_{AC_1}$, $3_B$, $3_{C_1}$\\
\hline 7& \small{all other partitions} & (0,0,0,0)  & \mbox{\small{all payoffs are} } $ 0 $\\
\hline
\end{tabular}
\end{center}
\label{dinner}
\end{table}%
The game runs as follows. Players simultaneously announce    individually chosen  coalition structures, what makes a point in a set of strategies of all players. Then a final coalition structure  is formed according to the rule above, and payoffs are assigned for the formed coalition structure. 

The way the game happens  means that a set of all strategies is divided by the  rule of coalition structure formation into domains: every domain corresponds to exactly one coalition structure. Or in other words, every  point in the set of strategies of the game  is assigned:  exactly one final coalition structure and payoffs  for every player. A player may have different payoffs in different   coalition  structures. Clearly  a final coalition structure  may not coincide with someone{'}s  initial choice.

Final payoffs have the following interpretation. Players A and B would always like  to be   together. Being rational they would choose a coalition structure with the highest payoff for them, i.e. the strategy 1. By the same reason the first best choices of $C_1$ and $C_2$  would  be   to choose  coalition structures  with  $A$. But A will never choose to be with either of them. The unavailability of the first best makes  $C_1$ and $C_2$ to choose option 2.\footnote{In sociology this behavior is referred as  a cooperation:    players $C_1$ and $C_2$ group together against other   options  when they are not together and have lower payoffs, player A will never choose them. This problem will be addressed in the accompanying paper.}
 By doing so they do not  disallow  a coalition $\{A,B \}$, but    reduce payoffs for A and   B. And players A and B cannot  prevent this (or to insure against). 

On the other hand, if players A and B choose strategy 2 they will  obtain coalition $\{A,B \}$ in any case, but in a different final coalition structure. In terms of mixed strategies this means that an equilibrium mixed strategy for $A$ and for $B$  is a whole probability space over two points, two  coalition structures 1 and 2.  For simplicity payoffs in coalition structures with coalition sizes three and four are much smaller.

From the forth column we can see that the corresponding cooperative game has an empty core. %, strong Nash equilibrium can not be applied here also, as well as coalition-proof equilibrium. 

The constructed game has a unique equilibrium. In terms of individual payoffs it is characterized as the second-best efficient for everybody. Equilibrium coalition structure contains   two coalitions. 
This equilibrium is different from  the strong Nash and a coalition proof equilibria. Both these concepts assume deviation of coalitions, but here definition of  a coalition is not enough to identify a final result of the game.  Additional problems appear with cores   - what is a direction of a deviation, may players deviate simultaneously from different coalitions, do deviators from different coalitions interact with each other?  Within cooperative game theory all these questions are left without well-defined answers. The proposed formal game in the paper  targets to  provide  coherent answers for these questions in terms of non-cooperative game theory.  And the result  of the game allows to study  inter / intra coalition externalities in a consistent non-cooperative way.

The  equilibrium in the corporate dinner game  does not   need super-additivity of payoffs.    %The game is based only on the  individual unbounded rationality of every  player.
There are  also differences from partition approach (Yi, 1999), as, for example, there is no initial allocation of players over coalitions.

% \end{document}
%\subsection{Reference to existing literature}
%A traditional way to study a coalition formation is through  a cooperative game theory.  However Maskin  (2011)    pointed out that   %"I believe that there are three (related)
%reasons for the historyorical 
%shift from cooperative to noncooperative theory: (a) most cooperative theory ignores externalities, the possibility that a coalition can be affected by the actions of those not in the coalition; (b) it assumes that a Pareto efficient outcome will be reached; and (c) it supposes that the grand coalition (the coalition of all players) will form."  
%Further: "These 
%"features of cooperative theory are problematic because most applications of game theory to economics involve settings in which externalities are important, Pareto inefficiency arises, and the grand coalition does not form. "  %This paper follows only  non-cooperative game theory approach.

\section{Formal setup  of the  model}
% Further in the paper we use the following terminology. Coalition structure (for short, a partition) is a collection of disjoint subsets  from  a set of players such that   a union of  all elements of a partition makes the original set.
%Coalition  is an element of a partition.

Nash (1950, 1951) suggested a non-cooperative game which consists of a set of players $N$, with a general element $i$, sets of  individual finite strategies $S_i,  i \in N$, and payoffs, defined as a mapping from a set of all strategies into  payoff    profiles of all players, $\Big(U_i(s)\Big)_{i \in N}$, such that    $S=\times_{i \in N}S_i \mapsto \Big(U_i(s)\Big)_{i \in N} \subset \mathbb{R}^{\#N}$, where $ \Big(U_i(s)\Big)_{i \in N} < \infty, \forall s \in S$.

The suggested game modifies the mapping by preliminary portioning $S$ into coalition structure specific domains and assigning payoffs for every point in these domains.  The division of $S$ is done with a coalition structure formation mechanism, defined  further.
  
Let there is a set of agents $N$, with a general element $i$, a size of  $N$  is  $\#N$, a  finite integer, $2 \le \#N < \infty$.

Every game has  a parameter $K$, $K \in \{1, \ldots, \#N\}$. This parameter has two interpretations. Let  for  $N$ agents there is a coalition with a maximum size  $K$. Then    no more than $K$  agents are required to dissolve it. %any  coalition of the size no more than $K$.   
The reverse is also true: we need no more than $K$ agents to form this maximum size  $K$ coalition.  
Closeness of   construction  of the object under investigation requires  these two simultaneous interpretations for   $K$ be equal.%: a maximum coalition size in  a  coalition structure  and a maximum number of deviators in the same coalition structure. %It is important to note that if we change value of $K$ we change the game.

Every value of $K$  from the set $\{1, \ldots, \#N \}$ induces a family of coalition structures  (or a family of partitions)   $\mathcal{P}(K)$ over the set of all players $N$: 
$$\mathcal{P}(K) = \{P=\{g_j \colon g_j \subset N; \#g \le K; \cup_j g_j = N;  \forall j_1 \neq j_2 \Rightarrow  {g}_{j_1} \cap {g}_{j_2} = \emptyset \}. $$

Every  coalition $g_j$ in a partition $P$ from $\mathcal{P}(K)$  has a size (a number of members) no bigger than $K$, but can be  less. The condition $j_1 \neq j_2 \Rightarrow  {g}_{j_1} \cap {g}_{j_2} = \emptyset $ is interpreted as  that an agent  can  participate only in one coalition. % from a partition. 

If we increase value of the parameter $K$ by one, then we need    to add  partitions from $\mathcal{P}(K+1) \setminus \mathcal{P}(K) $. This makes families of partitions for different $K$   be nested: % a set of eligible partitions,
  $\mathcal{P}(K=1) \subset \ldots \subset \mathcal{P}(K) \subset  \ldots  \subset \mathcal{P}(K=N)$.  The bigger is $K$, the more  coalition structures (or partitions) are  involved into consideration.  

For every partition $P$ an agent  $i$ has a finite strategy set $S_i (P)$.\footnote{Finite strategies are chosen as used in  Nash (1950),}   A set of strategies of agent $i$ for a family of coalition structures    $\mathcal{P}(K)$ is
$$
S_i(K) = \Big\{s_i(K) \colon s_i(K) \in \{ S_i(P) \colon P \in \mathcal{P}(K)\} \Big\}
$$ with a general element $s_i(K)$. 
For a given $K$ an agent chooses $s_i(K)$ from $S_i(K)$. \textit{ A choice of $s_i(K)$ means a choice of a desirable partition and an action for this partition}.\footnote{A desirable partition may not realize.  A coalition structure formation mechanism resolves conflicts of partition choices between players.} If we increase  the parameter $K$ by one, then we need to construct additional strategies only for the newly available coalition structures from  $\mathcal{P}(K+1) \setminus \mathcal{P}(K) $. This makes strategy sets for different $K$   be nested: % a set of eligible partitions,
  $S_i(K=1) \subset \ldots \subset S_i(K) \subset  \ldots  \subset S_i(K=N)$.

%For every family of coalition structures $\mathcal{P}(K)$ an agent $i$ has a  finite strategy set   $S_i(K) = \{S_i(P_i) \colon P_i \in \mathcal{P}(K)\}$, $P_i$ is a partition from $\mathcal{P}(K)$, which  $i$ may choose. %/network 
%actions of $i$  $P_i$ from $\mathcal{P}(K)$, % or strategies of $i$ for every choosable $P$ from  $\mathcal{P}$, 
%Individual  strategy sets for different $K$ are nested:
% $$S_i(K=1)  \subset \ldots \subset S_i(K) \subset  \ldots  \subset  S_i(K=N).$$ %Thus K creates a nested family of strategy sets.

 The set  of strategies for  a fixed $K$ is $S(K)=\times_{i \in N} S_i (K)$, a direct product of individual strategy sets of all players for  the given $K$. A choice of all players $s(K)=\Big(s_i(K), \ldots,s_N(K)\Big)$ is a point in   $S(K)$. For simplicity if there is no ambiguity we will write $s=\Big(s_1,\ldots,s_N\Big) \equiv s(K)=\Big(s_i(K), \ldots,s_N(K)\Big) $.
%Hence there are two ways to construct $S(K)$: in terms of initial individual strategies $S(K)=\times_{i}S_i (K)$, or in terms of realized partition $P$ specific strategies  $S(K)=\cup_{P} S(P) =\cup_{ P \in \mathcal{P}(K)} \cup_{i \in N}S_i(P)$. It is important that  a construction in terms of coalition structure specific strategies sets  may be not a direct product of sets, see an example in the next section.
It is clear that an increase in $K$ induces nested strategy sets:  $S(K=1) \subset \ldots  \subset S(K=N)$.

%Realization of a coalition structure may differ from someone{'}s choice. This may result in conflicts.
% All these cases the mechanism must anticipate and resolve at the stage of  construction of a game. The mechanism is embedded into the structure of the game.  

We have seen above that a mechanism to   resolve conflicts between individual choices is required. %It serves to   resolve conflicts between individual choices. 
   For every value  of $K$ from  the set  $  \{1,\ldots,\#N \}$ we define a coalition structure formation mechanism ( a mechanism  or a rule for short) $\mathcal{R}(K)$.  For  every  strategy profile $s=(s_1,\ldots,s_N) \in S(K)$ the mechanism assigns a final coalition structure $P$, $P \in \mathcal{P}(K)$, and transports the strategy profile $s$ into the partition $P$. Further   we will see that  this makes every $P$ to be a game, as $P$ will  have a set of players, a  non-trivial set of strategies and partition specific payoffs over it{'}s own strategy set.    %\footnote{and corresponding coalition structure strategy sets. This component will be illustrated in the next section example.} 
%We assume that for every $K$ there is a mapping which transports every chosen strategy profile $s$ from $S(K)$ into a strategy profile of  some final coalition structure $P$ from $\mathcal{P}(K)$.% and a corresponding strategy set of all players: 
%$\mathcal{R}(K) \colon \times_{i \in N}S_i(K) \mapsto \Big\{\cup_{i \in N} S_i(P) \colon P \in \mathcal{P}(K)\Big\}$. 

\begin{definition} For every $K$  \textbf{a coalition structure formation mechanism}  $\mathcal{R}(K)$  is a  set of mappings such that:
 \begin{enumerate}
\item A domain of $\mathcal{R}(K)$  is a set of all strategy profiles  of  $S(K)$. 
\item A range of $\mathcal{R}(K)$   is a finite  number of subsets  $S(P) \subset S(K)$, $P \in \mathcal{P}(K)$. $S(P)$ is a coalition structure specific strategy set .% , such that  $s \in S(P)$, where $S(P)$ is a coalition structure specific strategy set, $P \in \mathcal{P}(K)$.
\item Number of  the subsets $S(P)$ is no more than  a cardinality of $\mathcal{P}(K)$.%, such that $S(K)=\cup_{P \in \mathcal{P}(K)} S(P)$, where $S(P)$ is a coalition structure specific strategy set. 
\item   $\mathcal{R}(K)$ divides $S(K)$ into  coalition structure specific strategy sets, $ S(K) = \cup_{ P \in \mathcal{P}(K)}S(P)$. 
\item  Two different coalition structures, $\bar{P}$ and $\tilde{P}$, $ \bar{P} \neq \tilde{P}$,  have  different coalition structure   strategy sets    $  S( \bar{P}) \cap S(\tilde{P}) = \emptyset$.
\end{enumerate}
Formally the same:
 \begin{multline*}\mathcal{R}(K) \colon 
S(K)=\times_{i \in N}S_i(K)  \mapsto   \colon
 \begin{cases} 
 \forall s=(s_1,\ldots,s_N)  \in S(K)   \, \,   \exists P \in \mathcal{P}(K) \colon
 s \in S(P),    \\
 S(K) = \cup_{ P \in \mathcal{P}(K)}S(P), \\
  \forall \bar{P},\tilde{P} \in  \mathcal{P}(K),  \bar{P} \neq \tilde{P} \Rightarrow S( \bar{P}) \cap S(\tilde{P}) = \emptyset.
 \end{cases} .
\end{multline*}  
\end{definition}

Hence there are two ways to construct $S(K)$: in terms of initial individual strategies $S(K)=\times_{i \in N}S_i (K)$, and  in terms of   realized partition   strategies  $S(K) =\cup_{ P \in \mathcal{P}(K)} S(P)$. Representation of $S(K)$  in terms of coalition structure specific strategy  sets  may not be a direct product of sets, see an example in the next section. 
 We assume that a mechanism is given from outside. %, thus isolating questions of social welfare maximization from construction of a game.

If $K$ increases we need to add only  a mechanism for strategy sets  from   $S(K+1) \setminus S(K)$.   This supports  consistency of  coalition structure formation mechanisms for different $K$. 
The  family of mechanisms is nested:
: %not necessarily  what every $i$ chooses. 
 $\mathcal{R}(K=1)  \subset \ldots \subset \mathcal{R}(K) \subset  \ldots  \subset \mathcal{R}(K=N)$.  

 Payoffs in the game are defined as  state-contingent payoffs (or payoffs of Arrow-Debreu securities) in finance.  For every  coalition  structure  $P$   player $i$ has a payoff function $U_i(P) \colon S(P) \rightarrow \mathbb{R}_+$, such that  the set $U_i(P)$ is  bounded, $U_i(P)  < \infty$. Payoffs are considered as von  Neumann-Morgenstern utilities. All  payoffs  of $i$ for a game with no more than $K$ deviators make  a family:
  $\mathcal{U}_i(K) = \Big\{U_i(P) \colon   P \in \mathcal{P}(K) \Big\} $. Every coalition structure   has it{'}s own  set of strategies and a corresponding set of payoffs. Thus every coalition structure is a non-cooperative game.  
  
 An increase in $K$ also increases the number of possible partitions and the set of strategies for every player. We need to add only payoffs  for the partitions in  $\mathcal{P}(K+1) \setminus \mathcal{P}(K)$  or for strategies in $S_i(K+1) \setminus S_i(K)$. 
Thus   we obtain a nested family of payoff functions:
 $$ \mathcal{U}_i(K=1) \subset  \ldots \subset  \mathcal{U}_i(K) \subset  \ldots \subset \mathcal{U}_i(K=N).$$

 We can easily see that this construction of payoffs allows to obtain both intra and inter coalition (or group) externalities, as payoffs are defined  directly over strategy profiles of all players and independently from allocation of players in coalition structures. 

\begin{definition}[ \textbf{ a simultaneous coalition structure formation game}]
A non-cooperative game   for coalition structure formation is  
$$\Gamma(K)=\Big\langle N, \Big\{K, \mathcal{P}(K),\mathcal{R}(K) \Big\}, \Big(S_i(K),\mathcal{U}_i (K)\Big)_{i \in N} \Big\rangle, $$  where
$\Big\{K, \mathcal{P}(K),\mathcal{R}(K) \Big\}$  - coalition structure formation mechanism ( a social norm, a social institute),
$\Big(S_i(K),\mathcal{U}_i (K)\Big)_{i \in N}$ - properties of  players in $N$, (individual strategies and payoffs), such that: 
 $$\times_{i \in N} S_i(K)\buildrel \mathcal{R}(K) \over \rightarrow \Big\{ S(P) \colon P \in \mathcal{P}(K) \Big\} \rightarrow \Big\{(\mathcal{U}_i(K))_{i \in N}\Big\}. $$
\end{definition}

 For a example, the corporate dinner game above was a game for $K=2$, where an equilibrium did not change with an increase in values  of $K$  from $K=2$ to $K=3, 4$. 
 
Novelty of the paper is an introduction of coalition structure formation mechanism, which  portions the set of all strategies into  non-cooperative partition-specific  games. If we omit the mechanism part of the game and eliminate restriction on  coalition sizes then  we obtain the traditional non-cooperative game of Nash:  $\times_{i \in N} S_i(K)  \rightarrow  \Big(U_i(s)\Big)_{i \in N}. $
Construction of a game $\Gamma(K)$ makes every partition $P$  be an individual game. %, while all  partition-specific games in $\Gamma(K)$ are  connected by domains of mixed strategies.

Another novelty of the paper is an introduction of nested games.
 
\begin{definition}[\textbf{family of games}]
A family of games is \textbf{ nested}  if :
$$
\Gamma=\Gamma(K=1) \subset \ldots \subset \Gamma(K) \subset \ldots \subset \Gamma(K=N).
$$
 \end{definition}
 Nested games appear as a result of parametrization of a game by a maximum coalition size   (or by a maximum number of deviators, what is equivalent).  All games have consistent nesting of   components, for the same  set of  players $N$.%Existence of the nested property for games follows from constructions of all games in the family.
  
Let $\Sigma_i (K)$ be a set of all mixed strategies  of $i$, i.e. a space of probability measures,  $\Sigma_i (K) = \Big\{\sigma_i (K) \colon \int_{S_i(K)} d\sigma_i(K)=1 \Big\}$, with a general element $\sigma_i (K)$, where an integral is Lebegue-Stiltjes integral. Sets of mixed strategies for all other players are defined in the standard way $\Sigma_{-i} (K) = \Big\{\Big(\sigma_j(K)\Big)_{j \neq i} \colon  \forall j \neq i \int_{\Big(S_j(K)\Big)_{ j \neq i}} \Big(d\sigma_j(K)\Big) =1 \Big\}$.
  
Expected utility can be defined in terms of strategies  the  players  choose or in terms of final partition-specific strategies. 
Expected utility of $i$ in terms of individual strategies is: $$
EU^{\Gamma(K)}_i (\sigma_i(K),\sigma_{-i}(K)) = \int_{S(K)=\times_{i \in N} S_i(K)} {U}_i(s_i,s_{-i}) d\sigma_i(K) d\sigma_{-i}(K) $$
{ or in terms of partition-specific strategies is }

$$
EU^{\Gamma(K)}_i (\sigma_i(K),\sigma_{-i}(K)) = 
\sum_{P \in \mathcal{P}
(K)} \int_{S(P)} {U}_i (P)(s_i,s_{-i}) d\sigma_i(K)d\sigma_{-i}(K).$$ 

Expected utilities are constructed in the standard way.% 

 \begin{definition}[\textbf{ an equilibrium in a game $\Gamma(K)$ }]
A mixed strategies profile  $\sigma^*(K)=(\sigma^*_i(K))_{i \in N}$ is  an equilibrium strategy profile for a game $\Gamma(K)$ if for every  
  $\sigma_i(K) \neq \sigma^*_i(K)$ 
the following inequality for every player $i$ from $N$ holds true:
 $$
 EU^{\Gamma(K)}_i \Big(\sigma^*_i(K),\sigma^*_{-i}(K)\Big) \ge EU^{\Gamma(K)}_i \Big(\sigma_i(K),\sigma^*_{-i}(K)\Big).
  $$
\end{definition}

Equilibrium in the game  $\Gamma (K)$ is defined in a standard way. It{'}s existence   is just an expansion of Nash theorem.
%Every   game $\Gamma (K)$    from the family of non-cooperative games with coalition structure formations  $\mathcal{G}$ has  an equilibrium, may be in mixed strategies. 
However this result  for non-cooperative games with coalition structure formation is different from  the results  of cooperative games, where an equilibrium may not exist, for example,  in coalition form games with empty cores. Another outcome of the model  is  that there is no need to introduce additional properties of games, like axioms on a system of payoffs, super-additivity, or weights. Equilibrium existence result  can be generalized for the whole family of games.
  
\begin{theorem}
The family of games $\mathcal{G}=\{ \Gamma(K), K=1,2, \dots,N\}$ has  an equilibrium  in mixed strategies, $\sigma^*(\mathcal{G})=(\sigma^*(K=1), \ldots,\sigma^*(K=N) )$, $\Big(\sigma^*(K)\Big)_{i \in N}$.
\end{theorem}
Technical side off the result is obvious. The theorem expands the classic Nash theorem. %A proof is based on uniform convexity, boundness and continuity of a product of  sets of  mixed strategies  $\Big(\Sigma_i(K) \Big)_{i \in N}$ and application of a fixed point theorem.
An equilibrium in the game can also be  characterized  by   equilibrium partitions. 
\begin{definition}[\textbf{equilibrium coalition structures or partitions}]
A set of  partitions $\{ P^*\}(K)$, $\{ P^*\}(K) \subset \mathcal{P}(K)$, of a game $\Gamma(K)$, is a set of equilibrium partitions,  if it is induced by an equilibrium strategy profile $\sigma^*(K)=(\sigma^*_i(K))_{ i \in N}$. % with  domains of the equilibrium strategies  \textit{only} in these partitions:  $$Dom\Big(\sigma^*(K)\Big) \subset \{S(P^*)  \colon \{ P^*\}(K)  \subset \mathcal{P}(K)\}.$$
\end{definition}
In the same way we can define equilibrium payoffs for the whole family of games. This consideration is important for construction of a non-cooperative stability criterium in the next paper.

In the example below we will see  that there   can be more than one equilibrium partition, and equilibrium partitions may change with an increase in the number of deviators. These issues are addressed   in the accompanying paper on applications of the suggested game.

\section{An example of a game of two players}
The example serves to demonstrate two points. The constructed game has the standard property of games:   an equilibrium and efficiency concepts may not coincide, what is different from strong Nash concept. The second is to demonstrate that the concept of cooperation requires additional investigation (done  in the next  paper). The example below  is based on a   generalization of  the Prisoner{'}s  Dilemma. 

There are 2 players. They can form two types of partitions. If $K=1$ then there is only one final partition, $\mathcal{P}(K=1) \equiv P_{separ}=\{ \{ 1\},\{ 2\}\}$. If $K=2$ there are two final partitions, which make a family of partitions  $ \mathcal{P}(K=2)=\Big\{\{ \{ 1\},\{ 2\}\}, \{ 1, 2\} \Big\}$. Further we will use the notation $P_{joint}=\{ 1, 2\}$. Clearly, partition structures $\mathcal{P}(K=1)$ and $\mathcal{P}(K=2)$ are nested.

Every player in every partition has two strategies: H(igh) and L(ow). Player $i$ for a game with $K=1$ has a strategy set $S_i(K=1)=\{L_{i,P_{separ}},H_{i,P_{separ} }\}$. Player $i$ in the game with $ K=2$ has the strategy set $S_i(K=2)=\{ L_{i,P_{separ}},H_{i,P_{separ} }, L_{i,P_{joint}},H_{i,P_{joint} } \}.$ Clearly,   strategy sets   $S_i(K=1)$ and  $S_i(K=2)$ are nested.     

  Set of  strategies   for the  game  with $K=1$ is $ S(K=1)=\{L_{1,P_{separ}},H_{1,P_{separ} }\} \times \{L_{2,P_{separ}},H_{2,P_{separ} }\}$. Payoffs are in corresponding top-left cells   of  Table \ref{default2}.  Every cell   contains a payoff profile and a \textit{final} coalition structure. The payoffs  in these cells are   identical to    the payoffs of the Prisoner{'}s Dilemma.

   Set of  strategies   for the  game with $K=2$ is $$S(2)=\{ L_{1,P_{separ}},H_{1,P_{separ} }, L_{1,P_{joint}},H_{1,P_{joint} } \} \times \{ L_{2,P_{separ}},H_{2,P_{separ} }, L_{2,P_{joint}},H_{2,P_{joint} } \}.$$  Strategy sets of games with $K=1$ and $K=2$ are nested, i.e.  $S(K=1)  \subset S(K=2)$.  Additional payoffs of the game for $K=1$ are replicated from the game with $K=1$, see Table \ref{default2}. 
   
   For $K=1$ there are four outcomes, for $K=2$ there are sixteen outcomes, twelve new in comparison to $K=1$. 
  
 A coalition structure formation mechanisms for different $K$ are $\mathcal{R}(K=1),\mathcal{R}(K=2)$. The coalition structure formation mechanism is: the  partition $P_{joint}=\{1,2 \}$ can be formed only from a unanimous agreement of both  players.  It is  clear that  for $K=1$  only $P_{separ}$ can be formed. The grand coalition, or a partition $P_{joint}$, can be formed only over the strategy set  $(L_{1,P_{joint}},H_{1,P_{joint} })\times (L_{2,P_{joint}},H_{2,P_{joint} })$. Otherwise the partition $P_{separ}$ is formed. Thus a set of strategies for the partition $P_{separ}$  is not a direct product of some sets. 
 It is also clear that an increase in $K$ results in nested mechanisms: $\mathcal{R}(K=1) \subset \mathcal{R}(K=2)$. 
 
From Table \ref{default2} we can see that the whole  strategy set of the game  is partitioned into \textit{coalition structure specific domain\textbf{s}}.  Every coalition structure ( or a partition) is  a non-cooperative game  with it{'}s own strategy set and payoff profiles. Final partition for a player may not coincide with an individual choice.   A set of strategies for the partition $P_{separ}$ is not a product of sets. Finally  the games for $K=1$ and $K=2$ are nested.

Consider the game  with a maximum coalition size $K=2$ as described by the Table \ref{default2}.  If both  players choose strategies only for the partition $P_{separ}$, then the game   is the standard Prisoner{'}s Dilemma game. However the same coalition structure can be  formed if some of the  players does not choose any strategy for $P_{joint}=\{1,2\}$. %In other words there is a mismatch in choices. 
Thus for the \textit{final }  partition $P_{separ}=\{ \{1 \},\{ 2\}\}$ there are three equilibria, and every equilibrium is inefficient. There are also three efficient and desirable outcomes in the same partition. 

The partition $P_{joint}=\{1,2 \}$ can be formed only if both players choose it. Within this partition there is one   inefficient equilibrium and one efficient  outcome.

Compare efficient payoff profiles for the partitions $P_{separ}$ and $P_{joint}$. They have   equal  payoff  profiles $(0;0)$, but  their strategy profiles belong to different  final partitions. In other words, from   observing  only  the payoff profile $(0;0)$ we can not make a conclusion, which coalition structure is formed either $P_{separ}$ or $P_{joint}$. Another interpretation is that cooperation may take place as a cooperation in terms of positive externalities and in terms of allocation of players in one coalition. We can see that they are different.
%But  we can not address any of these efficient  outcomes as cooperative, as they appear in different partitions: one in $P_{separ}=\{\{1\},\{2 \}\}$  and another in $P_{joint}=\{1,2 \}$.

This means there are two kinds of cooperation - a cooperation in payoffs and a cooperation in allocation of players over coalitions. The example demonstrates that they may not coincide. More on that  in the next   paper. 
 
Using the same game we can demonstrate appearance of intra- and inter-coalitions externalities. 
If  partition $P_{joint}=\{1,2 \}$  is formed, then an individual payoff of a player depends on a strategy of another in the \textit{ same} coalition ( presence of intra-coalition  or intra-group externality). % and there is a positive additional gain from being together, i.e. $\epsilon >0$, both players are extravertes 
If partition $P_{separ}=\{\{1\} ,\{2 \}\}$, is formed then an individual payoff of a player depends on a strategy of another player  in  \textit{ the  different } coalition (  presence of inter-coalition  or inter-group externality). 

Multiplicity of equilibria makes  both of these externalities co-exist in equilibria, but in different final coalition structures.
Thus the game is able to present both intra and inter-coalition externalities, what is impossible in cooperative game theory.  

 \begin{table}[htp]
\caption{{Payoff for the family of games  with unanimous  formation rules}.  { {Different partitions have payoff-equal efficient outcomes.   }}
} 
\begin{center}
\begin{tabular}{|c|c|c||c|c|} \hline
 & $L_{2,P_{separ}}$ & $H_{2,P_{separ}} $& $L_{2,P_{joint}} $& $H_{2,P_{joint}} $ \\  \hline
$L_{1,P_{separ}}$ &
 \begin{tabular}{c} \underline{(0;0)}  \\ $ \{\{1\},\{ 2\} \}$ \end{tabular}

&  \begin{tabular}{c}  (-5;3) \\ $\{\{1\},\{ 2\} \}$  \end{tabular}

& \begin{tabular}{c} \underline{(0;0)} \\  $  \{\{1\},\{ 2\} \}$ \end{tabular}
& \begin{tabular}{c}  (-5;3) \\ $  \{\{1\},\{ 2\} \}$\end{tabular}  \\ \hline
$H_{1,P_{separ}}$ 
&  \begin{tabular}{c}  (3;-5) \\ $ \{\{1\},\{ 2\} \}$ \end{tabular}
&  \begin{tabular}{c}  \textbf{(-2;-2)} \\ $  \{\{1\},\{ 2\} \}$ \end{tabular}
&  \begin{tabular}{c}  (3;-5) \\ $  \{\{1\},\{ 2\} \}$ \end{tabular}
& \begin{tabular}{c}  \textbf{ (-2;-2)} \\ $  \{\{1\},\{ 2\} \}$\end{tabular} \\ 
\hline \hline
$L_{1,P_{joint}}$ 
&  \begin{tabular}{c}  \underline{(0;0)} \\ $  \{\{1\},\{ 2\} \}$\end{tabular}
&  \begin{tabular}{c}  (-5;3) \\ $  \{\{1\},\{ 2\} \}$\end{tabular}
&  \begin{tabular}{c}  $  \underline{(0;0)} $ \\ $  \{1, 2 \}$\end{tabular}
&  \begin{tabular}{c}  $(-5;3)$  \\ $  \{1,2\}$  \end{tabular}\\ \hline
$H_{1,P_{joint}}$ 
&  \begin{tabular}{c}  (3;-5) \\ $  \{\{1\},\{ 2\} \}$ \end{tabular}
&  \begin{tabular}{c}  \textbf{ (-2;-2)} \\ $  \{\{1\},\{ 2\} \}$ \end{tabular}
& \begin{tabular}{c}   $(3;-5) $  \\ $  \{1,2\}$ \end{tabular}
&  \begin{tabular}{c}  \textbf{ (-2;-2)} \\ $  \{1, 2\}$   \end{tabular}\\ \hline  
 \end{tabular}
\end{center}
\label{default2}
\end{table}%
  
{\tiny{\begin{table}[htp]
\caption{Payoff for two extrovert players who, obtain additional payoffs  $\epsilon$ being in one coalition $P_{joint}$, when it is realized. Uniqueness of  an equilibrium is reinstalled.}
\begin{center}
\begin{tabular}{|c|c|c||c|c|} \hline 
 & $L_{2,P_{separ}}$ & $H_{2,P_{separ}} $& $L_{2,P_{joint}} $& $H_{2,P_{joint}} $ \\  \hline
$L_{1,P_{separ}}$ &
 \begin{tabular}{c} \underline{(0;0)}  \\ $ \{\{1\},\{ 2\} \}$ \end{tabular}

&  \begin{tabular}{c}  (-5;3) \\ $\{\{1\},\{ 2\} \}$  \end{tabular}

& \begin{tabular}{c} \underline{(0;0)} \\  $  \{\{1\},\{ 2\} \}$ \end{tabular}
& \begin{tabular}{c}  (-5;3) \\ $  \{\{1\},\{ 2\} \}$\end{tabular}  \\ \hline
$H_{1,P_{separ}}$ 
&  \begin{tabular}{c}  (3;-5) \\ $ \{\{1\},\{ 2\} \}$ \end{tabular}
&  \begin{tabular}{c}  \textbf{(-2;-2)} \\ $  \{\{1\},\{ 2\} \}$ \end{tabular}
&  \begin{tabular}{c}  (3;-5) \\ $  \{\{1\},\{ 2\} \}$ \end{tabular}
& \begin{tabular}{c}  \textbf{ (-2;-2)} \\ $  \{\{1\},\{ 2\} \}$\end{tabular} \\ 
 \hline \hline
$L_{1,P_{joint}}$ & 
\begin{tabular}{c} \underline{(0;0)} \\ $ \{\{1\},\{ 2\} \}$\end{tabular}
& \begin{tabular}{c} (-5;3) \\ $  \{\{1\},\{ 2\} \}$ \end{tabular}
& \begin{tabular}{c} $  \underline{(0+\epsilon;0+\epsilon)} $ \\  $\{1, 2 \}$ \end{tabular}& \begin{tabular}{c} $(-5+\epsilon;3+\epsilon)$ \\  $\{1,2\}$ \end{tabular} \\ \hline
$H_{1,P_{joint}}$ 
& \begin{tabular}{c} (3;-5) \\ $ \{\{1\},\{ 2\} \}$\end{tabular} 
&  \begin{tabular}{c} (-2;-2)  \\ $  \{\{1\},\{ 2\} \}$ \end{tabular} & 
\begin{tabular}{c}
$(3+\epsilon;-5+\epsilon) $ \\ $ \{1,2\}$ \end{tabular} &  \begin{tabular}{c} \textbf{ (-2+$\epsilon$;-2+$\epsilon$)} \\ $ \{1,2 \}$  \end{tabular} \\ \hline  
 \end{tabular}
\end{center}
\label{default3}
\end{table}%
}}

We can reinstall uniqueness of an equilibrium, what is done in Table \ref{default3}.  If  both players are extroverts and prefer be together,\footnote{what is equivalent to preferences over coalition structures} then  every individual payoff increases by $\epsilon>0$, if the grand coalition is realized. This means that a  change of   the game from $\Gamma(1)$  to $\Gamma(2)$ changes the equilibrium in terms of both   strategies and the   partitions.

 If both players are introverts,  $\epsilon <0$, then the expansion of the game will not change initial equilibrium in terms of both   strategies and the   partition. This case is not  presented  here and is left for future. %as it directly leads to trembling hand equiliHowever in both cases for extroverts and introverts there are no changes in   equilibrium payoff profiles. These  issues are discussed in  the  accompanying paper. 

\section{Discussion}

%Game in Nash (1950,51):  $\times_{i \in N} S_i(K) \rightarrow \{(\mathcal{U}_i(K))_{i \in N}\}. $

%Family of games in this  paper  $\times_{i \in N} S_i(K)\buildrel \mathcal{R}(K) \over \rightarrow \{ S(P) \colon P \in \mathcal{P}(K) \} \rightarrow \{(\mathcal{U}_i(K))_{i \in N}\}, K = 1,N. $

Insufficiency of cooperative game theory to study  coalitions and coalition structures was earlier reported by many authors. Maskin (2011) wrote that \textquotedblleft  features of cooperative theory are problematic because most applications of game theory to economics involve settings in which externalities are important, Pareto inefficiency arises, and the grand coalition does not form\textquotedblright.  
 Myerson (p.370, 1991) noted that  \textquotedblleft  we need some model of cooperative behavior that does not abandon the individual decision-theoretic foundations of game theory\textquotedblright. Thus there is a  demand for   a specially designed  non-cooperative game to study coalition structures formation along  with an adequate equilibrium concept for this game.  
 
 There is a voluminous literature on the topic, a list of authors is far from complete:  Aumann, Hart, Holt,  Maschler,  Maskin, Myerson, Peleg,  Roth, Serrano, Shapley, Schmeidler, Weber, Winter, Wooders and many others.

A popular approach to use   a \textquotedblleft threat\textquotedblright   as a basic concept for  coalition   formation analysis
 was suggested by Nash (1953).    
  %One way  to structure is to construct threats (in a terminology of Nash, 1953).
   Consider a  strategy profile from  a subset of   players. Let this profile be  a threat to someone, beyond this subset.  
The threatening player\textbf{s}  may  produce externalities   for each other (and negative  externalities not excluding!). How credible  could  be such threat? At the same time  there may be some other player beyond the subset of players who may obtain a bonanza from this threat.  But  this beneficiant may  not join the group due to expected intra-group negative externalities  for members  or from members of this group.  Thus a  concept of a threat can not serve as an elementary concept.

 The justification of a chosen tool, a  non-cooperative game, comes from   Maskin (2011) and a     remark of 
 Serrano (2014), %on a preference of   non-cooperative to cooperative game theory:  "the axiomatic route find difficulties identifying solutions", and 
 that for studying coalition formation \textquotedblleft it  may be worth  to use strategic-form   games, as proposed in the Nash program\textquotedblright.

The difference of the research agenda   in this paper  from  the  Nash program (Serrano 2004) is   studying   non-cooperative formation  of  coalition structures, but  not only  formation of  one coalition. The best analogy for the difference   is the  difference between partial and general strategic equilibrium analysis in economics. The former isolates a market  ignoring cross-market interactions, the latter explicitly studies cross market interactions.

The constructed finite non-cooperative game allows to  study  what can be a cooperative behavior, when  the individuals
 %The notion that groups of individuals will act to achieve their common or group interests, far from being a logical implication of the assumption that 
 \textquotedblleft    rationally further their individual interests" ( Olson, 1971). %is in fact inconsistent with that assumption... If the members of a large group rationally seek to maximize their personal welfare, they will not act to advance their common or group objectives unless there is coercion to force them to do so, or unless some separate incentive, distinct from the achievement of the common or group interest, is offered to the members of the group individually on the condition that they help bear the costs or burdens involved in the achievement of the group objectives." 
 
 Nash (1950, 1951) suggested to construct a non-cooperative game as a mapping of a set of strategies into a profile of payoffs,   $\times_{i \in N} S_i  \rightarrow (U_i)_{i \in N}. $ 

This paper has two contributions in comparison to his paper:  construction of a non-cooperative game with an  embedding  coalition structure formation mechanism, and   parametrization of all constructed games by a number of deviators:   $  \times_{i \in N} S_i(K)\buildrel \mathcal{R}(K) \over \rightarrow \{ S(P) \colon P \in \mathcal{P}(K) \} \rightarrow \{(\mathcal{U}_i(K))_{i \in N}\}, $ where $K \in \{ 1,\ldots,\#N \} $. The game suggested by Nash becomes a partial case for these games.

  Every game in a family   has an equilibrium, may be in mixed strategies. This differs from results of cooperative game theory, where  games may have no equilibrium, like  games with empty cores, etc. 
 
 The introduced  equilibrium concept differs from the strong  Nash, coalition-proof and $k$ equilibrium concepts. The differences  are:  an explicit allocation of payoffs and a combined presence of intra- and inter- coalition (or group) externalities (the list  of differences is not complete). Differences from the core approach of Aumann (1960) are clear:  a presence of externalities, no restrictions that  only one group deviates, no restrictions on the direction of a deviation (inside or outside), and a  construction of individual payoffs from a strategy profile of all players. There is no need to assume transferable utilities for players.
The approach allows to  study    coalition structures, which differ  from the grand coalition  as in   Shapley value.  
 Finally the introduced concepts enables to offer a non-cooperative  necessary stability criterion based only on an equilibrium of a game. This is  impossible for any other equilibrium concept.

The suggested  approach is different  in a role for a central planner offered by Nash, who \textquotedblleft  argued that cooperative actions are the result of some process of bargaining" Myerson (p.370, 1991).  
%Self-interest cooperative actions  in the model appear  from a strategic non-cooperative behavior of every player within  a predefined mechanism design framework. 
The central planner offers a predefined coalition structure formation mechanism, that includes a maximum number of deviators,  family of eligible partitions and a family of rules  to construct these partitions from individual strategies of players.

 The  accompanying papers will demonstrate application of the suggested model for studying stochastic Bayesian games, cooperation, self-enforcement properties of an equilibrium  (Aumann, 1990), non-cooperative criterion of partition stability,  focal points, application of the same mechanism to study network games and repeated games.

\end{document}